\documentclass[11pt]{amsart}
\usepackage{}
\usepackage{amssymb}
\theoremstyle{plain}
\newtheorem{Thm}{Theorem}

\begin{document}

\title[gradient estimates for a nonlinear heat equation]
{Gradient estimates for a simple nonlinear heat equation on
manifolds}

\author{Li Ma}

\address{Department of mathematics \\
Henan Normal university \\
Xinxiang, 453007 \\
China}

\email{nuslma@gmail.com}

\thanks{The research is partially supported by the National Natural Science
Foundation of China 10631020 and SRFDP 20090002110019}

\begin{abstract}
In this paper, we study the gradient estimate for positive
solutions to the following nonlinear heat equation problem
$$
u_t-\Delta u=au\log u+Vu, \ \  u>0
$$
on  the compact Riemannian manifold $(M,g)$ of dimension $n$ and
with non-negative Ricci curvature. Here $a\leq 0$ is a constant,
$V$ is a smooth function on $M$ with $-\Delta V\leq A$ for some
positive constant $A$. This heat equation is a basic evolution
equation and it can be considered as the negative gradient heat
flow to $W$-functional (introduced by G.Perelman), which is the
Log-Sobolev inequalities on the Riemannian manifold and $V$
corresponds to the scalar curvature.


\end{abstract}

 \maketitle

\section{Introduction}

The study of the heat equation
\begin{equation}\label{ht}
u_t-\Delta u=0, \ \  u>0
\end{equation}
on the compact Riemannian manifold $(M,g)$ of dimension $n$ plays
an important role in the field of Partial differential equations,
Geometric analysis, and Differential Geometry (\cite{Chow}).
Rather than mentioning a lot references, we prefer to give a
comment about our appendix studying the quantity introduced by
G.Perelman for the gradient estimate of the fundamental solution
to the conjugate heat equation associated with the Ricci flow.
From it, one can see that the choices of differential Harnack
quantities should depend on the geometric background of eveolution
equations. Because Perelman's gradient estimate for heat kernel
plays an important role in Ricci flow, people are motivated to
find extensions of his result to various problems. In \cite{Ni} L.
Ni made a clever observation that the same gradient estimate for
the fundamental solution to the heat equation on a compact
Riemannian manifold with fixed metric is also true provided the
Ricci curvature is non-negative. Recently attention is focused to
a nature question, which asks if such a gradient estimate is still
true for all positive solutions to heat equation or conjugate heat
equation associated with the Ricci flow. Cao-Hamilton\cite{CH} and
Kuang-Zhang \cite{ KZ} then showed independently that the
Perelman's gradient estimate is true for the conjugate heat
equation associated with the Ricci flow. In our appendix, we
notice that the Perelman type gradient estimate is an easy
consequence of Li-Yau gradient estimate for any positive solution
to the heat equation on a compact Riemannian manifold with
non-negative Ricci curvature. Actually the above estimate can be
obtained easily from Li-Yau's gradient estimate (\cite{LY}),
namely,
$$
2t\Delta f\leq n.
$$
Just note that $$ tP\leq tP+t|\nabla f|^2=2t\Delta f\leq n,$$
where the constant $n$ is the best. Other method may not give such
a sharp constant. Following our works \cite{M} and \cite{M2}, we
continue to study the gradient estimate for positive solutions to
nonlinear evolution equations on manifolds.

 The main topic of this note is to study the following nonlinear heat
equation problem
\begin{equation}\label{n-ht1}
u_t-\Delta u=au\log u+Vu, \ \  u>0
\end{equation}
on  the compact Riemannian manifold $(M,g)$ of dimension $n$. Here
$a\leq 0$ is some constant and $V$ is a given smooth function on
$M$. This heat equation can be considered as the negative gradient
heat flow to $W$-functional \cite{P02}, which is the Log-Sobolev
inequalities on the Riemannian manifold and $V$ corresponds to the
scalar curvature.

Our result is below.

\begin{Thm}\label{thm2}
Assume that the compact Riemannian manifold $(M,g)$ has
non-negative Ricci curvature. Assume $-\Delta V\leq A$ for some
constant $A\geq 0$.
 Let $u>0$ be a positive smooth
solution to (\ref{n-ht1}). Let $f=\log u$. Then we have, for all
$t>0$,
$$
\Delta f-At-\frac{n}{2t}\leq 0,
$$
and in other words, $$ f_t-af+V+|\nabla f|^2\leq At+\frac{n}{2t}.
$$ The same result is also true for complete Riemannian manifold
provided the maximum principle is applicable.
\end{Thm}

The important part in our result is that we only assume the lower
bound of $\Delta V$ on $M$. However, our gradient estimate may not
be the best one because of the term $At$. In our paper \cite{M},
we propose the study the nonlinear heat equation (\ref{n-ht1}).
Then Y.Yang \cite{Y} finds a nice gradient estimate for
(\ref{n-ht1}). G.Huang \cite{H}, Chen and Chen \cite{CC}, and
others also find more interesting results for (\ref{n-ht1}). See
also \cite{M2} for related rsults.

We can extend the result above to compact Riemannian manifold with
smooth convex boundary.

\begin{Thm}\label{thm3}
Assume that the compact Riemannian manifold $(M,g)$  with smooth
convex boundary has non-negative Ricci curvature. Assume $-\Delta
V\leq A$ for some constant $A\geq 0$.
 Let $u>0$ be a positive smooth
solution to (\ref{n-ht1}) with Neumann boundary condition $u_\nu$
where $\nu$ is the outward unit normal to the boundary. Assume
that $V_\nu\leq 0$ on the boundary $\partial M$. Let $f=\log u$.
Then we have, for all $t>0$,
$$
\Delta f-At-\frac{n}{2t}\leq 0.
$$
\end{Thm}

From the geometric view-point, the boundary condition about $V$ is
not nature. We point out that it is interesting question to extend
the results above to the equation (\ref{n-ht1}) on  general
Riemannian manifolds.

\section{Proofs of Theorems \ref{thm2} and \ref{thm3}}

We shall follow Li-Yau's method \cite{LY}.

Consider the positive solution $u$ to the equation (\ref{n-ht1}).
Let $f=-\log u$. Then we have
$$
Lf:=f_t-\Delta f=af-V-|\nabla f|^2.
$$
Then we have
$$
L\Delta f=a\Delta f-\Delta V-\Delta |\nabla f|^2.
$$
By the Bochner formula and the non-negative Ricci curvature
assumption we know that
$$
\Delta |\nabla f|^2\geq \frac{2}{n}|\Delta f|^2+2g(\nabla f,\nabla
\Delta f).
$$

Then we have
$$
L\Delta f-a\Delta f\leq -\Delta V-\frac{2}{n}|\Delta
f|^2-2g(\nabla f,\nabla \Delta f).
$$

Recall that $-\Delta V\leq A$ for some constant $A\geq 0$. Hence
we have
$$
L\Delta f-a\Delta f\leq A-\frac{2}{n}|\Delta f|^2-2g(\nabla
f,\nabla \Delta f).
$$
Define
$$
Q=\Delta f-At-\frac{n}{2t}.
$$
We then have
$$
LQ-aQ +2g(\nabla f,\nabla Q)\leq a(At+\frac{n}{2t})+
\frac{2}{n}[(\frac{n}{2t})^2-|\Delta f|^2].
$$
Since $a\leq 0$, we have
$$
LQ-aQ +2g(\nabla f,\nabla Q)\leq
\frac{2}{n}[(\frac{n}{2t})^2-|Q+At+\frac{n}{2t}|^2].
$$
Note that
$$
[(\frac{n}{2t})^2-|Q+At+\frac{n}{2t}|^2]=-(Q+At)(Q+At+\frac{n}{t}).
$$

\textbf{Case 1.} At point $(x,t)$ where
$$
Q+At+\frac{n}{t}\leq 0,
$$
we have $Q+At\leq 0$ too. Then we have
$$
LQ-aQ +2g(\nabla f,\nabla Q)\leq 0.
$$

\textbf{Case 2.} At point $(x,t)$ where
$$
Q+At+\frac{n}{t}\geq 0,
$$
we have $$ -(Q+At)(Q+At+\frac{n}{t}) \leq -Q(Q+At+\frac{n}{t})
$$
and then
$$
LQ-aQ +2g(\nabla f,\nabla Q)\leq -Q(Q+At+\frac{n}{t}).
$$

Define the function $B=0$ at the point where $Q+At+\frac{n}{t}\leq
0$ and $B=-Q(Q+At+\frac{n}{t})$ at the point where
$Q+At+\frac{n}{t}\geq 0$. Then using the maximum principle we
obtain that
$$
Q:=\Delta f-At-\frac{n}{2t}\leq 0.
$$
This completes the proof of Theorem \ref{thm2}.

The proof of Theorem \ref{thm3} is similar. We need only to
exclude the possibility of the maximum point of $Q$ at boundary
points. If the maximum occurs at the boundary point $(x_0,t_0)$,
then by the strong maximum principle we have
$$
Q_\nu>0
$$
at this point. Note that $$\Delta f= f_t-af+V+|\nabla f|^2.$$ So,
at $(x_0,t_0)$,
$$
Q_\nu=(\Delta f)_\nu=V_\nu+(|\nabla f|^2)_\nu\leq -2II(\nabla
f,\nabla f)\leq 0.
$$
A contradiction. This proves Theorem \ref{thm3}.

\section{Appendix}

To compare the Perelman and Li-Yau gradient estimate for heat
equation, we follow Perelman's approach to study the equation
(\ref{ht}). Precisely, we prove the following.

\begin{Thm}\label{thm1}
Assume that the compact Riemannian manifold $(M,g)$ has
non-negative Ricci curvature. Let $u>0$ be a positive smooth
solution to (\ref{ht}). Let $f=\log u$. Then we have, for all
$t>0$,
$$
P:=2\Delta f-|\nabla f|^2\leq 2n/t.
$$
\end{Thm}

We now prove Theorem \ref{thm1}. Recall that $f=-\log u$. Then
$$
f_j=-u_j/u, \ \ \Delta f=-\Delta u/u+|\nabla f|^2.
$$
Then we have
$$
(\partial_t-\Delta)f=-|\nabla f|^2.
$$

Let $L=\partial_t-\Delta$.
 Set $P=2y-z$ with $y=\Delta f$ and
$z=|\nabla f|^2$.

Note that
$$
Lf=-z.
$$

Then $Ly=-\Delta z$. Compute
$$
Lz=-\Delta z+2g(\nabla f,\nabla f_t)
$$
$$
=-\Delta z+2g(\nabla f,\nabla (y-z)).
$$

Using the Bochner formula
$$
\Delta z=2|D^2f|^2+2g(\nabla f,\nabla y)+2Rc(\nabla f,\nabla f),
$$
 we get that
$$
LP=-2g(\nabla f,\nabla P)-2|D^2f|^2-2Rc(\nabla f,\nabla f).
$$

Using the Cauchy-Schwartz inequality, we obtain that
$$
|D^2f|^2\geq \frac{1}{4n}(P+z)^2.
$$

Hence using $Rc\geq 0$, we have
$$
LP+2g(\nabla f,\nabla P)\leq -\frac{1}{2n}(P+z)^2
$$
and
$$
L(P-2nt^{-1})+2g(\nabla f,\nabla (P-2nt^{-1}))\leq
-\frac{1}{2n}[(P+z)^2-(2n/t)^2].
$$

Note that
$$
I:=[(P+z)^2-(2n/t)^2]=[P+z-(2n/t)][P+z+(2n/t)].
$$

Then we have
$$
I=[P+z-(2n/t)]^2+\frac{4n}{t}[P+z-(2n/t)]\geq
\frac{4n}{t}[P+z-(2n/t)].
$$
Note that $z=|\nabla f|^2\geq 0$. We then have
$$
I\geq \frac{2}{t}[P-(2n/t)].
$$
Then we have
$$
L(P-2nt^{-1})+2g(\nabla f,\nabla (P-2nt^{-1}))\leq
-\frac{1}{2n}\frac{2}{t}[P-(2n/t)].
$$
Using the maximum principle on $M\times [\delta, T)$ for
$\delta>0$ small, we conclude that
$$
P-(2n/t)\leq 0.
$$
This completes the proof of Theorem \ref{thm1}.

We remark that one may also use the argument in \cite{KZ} in the
following way.

 \textbf{Case 1.} Assume that
$$
P+z+(2n/t)\leq 0.
$$

Then
$$
P+z-(2n/t)\leq 0
$$
too and $I\geq 0$.

\textbf{Case 2.}

 Assume that
$$
P+z+(2n/t)\geq 0.
$$

Then
\begin{align*}
I&=z[P+z+(2n/t)]+[P-(2n/t)][P+z+(2n/t)]\\
&\geq [P-(2n/t)][P+z+(2n/t)].
\end{align*}

Define a new function $B=0$ on the set where $P+z+(2n/t)\leq 0$
and $B=P+z+(2n/t)$ on the set where $P+z+(2n/t)\geq 0$. Then we
have
$$L(P-2nt^{-1})+2g(\nabla f,\nabla (P-2nt^{-1}))\leq
-\frac{1}{2n}B[P-(2n/t)].
$$
Using the maximum principle again to get the conclusion as before.


\begin{thebibliography}{20}

\bibitem{CH}
X.Cao, R.Hamilton, \emph{Differential Hrnack estimate for
time-dependent heat equations with potential}, Geom. Funct. Anal.
19 (2009), no. 4, 989-1000

\bibitem{CZ} X.Cao, Q.Zhang,
\emph{The Conjugate heat equation and ancient solutions of the
Ricci flow}, Arxiv:1006.0540v1

\bibitem{CC}
Chen, Li; Chen, Wenyi, Gradient estimates for a nonlinear
parabolic equation on complete non-compact Riemannian manifolds.
Ann. Global Anal. Geom. 35 (2009), no. 4, 397-404.

\bibitem{Chow} Bennett Chow, Sun-Chin Chu, David Glickenstein,
Christine Guenther, James Isenberg, Tom Ivey, Dan Knopf, Peng Lu,
Feng Luo, and Lei Ni. \emph{The Ricci flow: techniques and
applications}. Part III. Mathematical Surveys and Monographs.
American Mathematical Society, Providence, RI. Geometric- Analysis
aspects, 2010.

\bibitem{H}
Huang, Guangyue; Ma, Bingqing, \emph{Gradient estimates for a
nonlinear parabolic equation on Riemannian manifolds}. Arch. Math.
(Basel) 94 (2010), no. 3, 265-275.

\bibitem{KZ}
S.Kuang, Qi. Zhang, \emph{A gradient estimate for all positive
solutions of the conjugate heat equation under Ricci flow}, J.
Funct. Anal. 255 (2008), no. 4, 1008-1023.

\bibitem{LY}
P. Li, S.T. Yau, \emph {On the parabolic kernel of the
Schr\"{o}inger operator}, Acta Math. 156 (1986) 153-201.

\bibitem{M}
Ma, Li, \emph{Gradient estimates for a simple elliptic equation on
complete non-compact Riemannian manifolds}. J. Funct. Anal. 241
(2006), no. 1, 374-382.

\bibitem{M2}
Ma, Li; Zhao, Lin; Song, Xianfa, \emph{Gradient estimate for the
degenerate parabolic equation $u_t=\Delta F(u)+H(u)$ on
manifolds}. J. Differential Equations 244 (2008), no. 5,
1157-1177.

\bibitem{Ni}
L.Ni, \emph{The entropy formula for linear Heat equation}, J.
Geom. Anal. 14(2004)87-100. Addenda to "he entropy formula for
linear Heat equation", J. Geom. Anal. 14(2004)369-374

\bibitem{P02} Grisha Perelman,
\emph{The entropy formula for the Ricci flow and its geometric
applications}, http://arxiv.org/abs/math/0211159v1


\bibitem{Y}
Yang, Yunyan, \emph{Gradient estimates for a nonlinear parabolic
equation on Riemannian manifolds}. Proc. Amer. Math. Soc. 136
(2008), no. 11, 4095-4102.












\end{thebibliography}
\end{document}